\documentclass[11pt]{article} 
\usepackage{amsmath,amssymb,amsthm} 
\usepackage{fancybox}  
\usepackage{graphicx} 
\usepackage{color}
\usepackage{colortbl}

\allowdisplaybreaks
                                                          
\begin{document}

\def\fl#1{\left\lfloor#1\right\rfloor}
\def\cl#1{\left\lceil#1\right\rceil}
\def\ang#1{\left\langle#1\right\rangle}
\def\stf#1#2{\left[#1\atop#2\right]} 
\def\sts#1#2{\left\{#1\atop#2\right\}}
\def\eul#1#2{\left\langle#1\atop#2\right\rangle}
\def\N{\mathbb N}
\def\Z{\mathbb Z}
\def\R{\mathbb R}
\def\C{\mathbb C}
\newcommand{\ctext}[1]{\raise0.2ex\hbox{\textcircled{\scriptsize{#1}}}}

\newtheorem{theorem}{Theorem}
\newtheorem{Prop}{Proposition}
\newtheorem{Cor}{Corollary}
\newtheorem{Lem}{Lemma}
\newtheorem{Def}{Definition}
\newtheorem{Rem}{Remark}

\title{The $p$-numerical semigroup of the triple of arithmetic progressions
}

\author{
Takao Komatsu 
\\
\small Department of Mathematical Sciences, School of Science\\[-0.8ex]
\small Zhejiang Sci-Tech University\\[-0.8ex]
\small Hangzhou 310018 China\\[-0.8ex]
\small \texttt{komatsu@zstu.edu.cn}\\\\
Haotian Ying
\\ 
\small Department of Mathematical Sciences, School of Science\\[-0.8ex]
\small Zhejiang Sci-Tech University\\[-0.8ex]
\small Hangzhou 310018 China\\[-0.8ex]
\small \texttt
{tomyinght@gmail.com} 
}

\date{
\small MR Subject Classifications: Primary 11D07; Secondary 05A15, 05A17, 05A19, 11B68, 11D04, 11P81, 20M14 
}

\maketitle
 
\begin{abstract} 
For given positive integers $a_1,a_2,\dots,a_k$ with $\gcd(a_1,a_2,\dots,a_k)=1$, the denumerant $d(n)=d(n;a_1,a_2,\dots,a_k)$ is the number of nonnegative solutions $(x_1,x_2,\dots,x_k)$ of the linear equation $a_1 x_1+a_2 x_2+\dots+a_k x_k=n$ for a positive integer $n$. 
For a given nonnegative integer $p$, let $S_p=S_p(a_1,a_2,\dots,a_k)$ be the set of all nonnegative integers $n$'s such that $d(n)>p$. In this paper, we are interested in the $p$-Frobenius number, which is the maximum of the set of gaps $\mathbb N_0\backslash S_p$. Here $\mathbb N_0$ denotes the set of nonnegative integers.  
When $p=0$, $S=S_0$ is the original numerical semigroup, and the $0$-Frobenius number is the original Frobenius number. 
The explicit formula for two variables is known not only for $p=0$ but also for $p>0$, but when there are three or more variables, it is difficult even in the special case of $p=0$. For $p>0$, it is not only more difficult, but no explicit formula had been found. In this paper, explicit formulas of the $p$-Frobenius number and related values are given for the triple of arithmetic progressions. The main tool is to determine the elements of the $p$-Ap\'ery set.     
\end{abstract}

\section{Introduction}  

For integer $k\ge 2$, consider a set of positive integers $A=\{a_1,\dots,a_k\}$. Denote by $d(n)=d(n;a_1,\dots,a_k)$ the number of nonnegative solutions $(x_1,x_2,\dots,x_k)$ of the linear equation $a_1 x_1+a_2 x_2+\dots+a_k x_k=n$ for a positive integer $n$. Recently, the concept of $p$-numerical semigroups is introduced together with their symmetric characteristics \cite{KY24}. $d(n)$ is often called the {\it denumerant}. For a given nonnegative integer $p$, let $S_p=S_p(A)$ be the set of all nonnegative integers $n$'s such that $d(n)>p$. 
For the set of nonnegative integers $\mathbb N_0$, the set $\mathbb N_0\backslash S_p$ is finite if and only if $\gcd(a_1,\dots,a_k)=1$. Then there exists the largest element $g_p(A)$ in $\mathbb N_0\backslash S_p$, which is called the {\it $p$-Frobenius number}, and each element is called the {\it gap}. For the so-called {\it $p$-numerical semigroup} $S_p(A)$, the cardinality $n_p(A)$ and the sum $s_p(A)$ of $\mathbb N_0\backslash S_p$ are called {\it $p$-Sylvester number} (or the {\it $p$-genus}) and the {\it $p$-Sylvester sum}, respectively.  
When $p=0$, $S=S_0(A)$ is the original numerical semigroup, and the $0$-Frobenius number $g_0(A)$is the original Frobenius number $g(A)$. Finding the Frobenius and related values is the well-known linear Diophantine problem, posed by Sylvester \cite{sy1884}, but known as the Frobenius problem, is the problem to determine the Frobenius number $g(A)$. The Frobenius Problem has been also known as the {\it Coin Exchange Problem} (or Postage Stamp Problem / Chicken McNugget Problem), which has a long history and is one of the problems that has attracted many people as well as experts. The genus $g(A)=g_0(A)$ is often fundamental in the study of algebraic curves and commutative algebra.  

For two variables $A=\{a,b\}$, it is shown that  
\begin{equation}
g(a,b)=(a-1)(b-1)-1\quad\hbox{and}\quad n(a,b)=\frac{(a-1)(b-1)}{2}\,. 
\label{eq:g-n}
\end{equation}   
(\cite{sy1882,sy1884}). 
An explicit expression of the Sylvester sum $s(A)=s_0(A)$ is given by Brown and Shiue \cite{bs93} for two variables $A=\{a,b\}$ as  
\begin{equation}
s(a,b)=\frac{1}{12}(a-1)(b-1)(2 a b-a-b-1)\,. 
\label{eq:s}
\end{equation}  
This result is extended in \cite{ro94} for the power sum of the set of gaps $s^{(\mu)}(A)$, defined by 
$$
s^{(\mu)}(A)=\sum_{n\in\mathbb N_0\backslash S_p(A)}n^\mu 
$$ 
in the case of $A=\{a,b\}$.  
However, for three or more variables, it very complicated to find a general explicit formula for Frobenius number, Sylvester number and Sylvester sum. Only for some special cases, explicit formulas have been found, including arithmetic, geometric-like, Fibonacci, Mersenne, and triangular (see, e.g., \cite{RR18} and references therein). The study of semigroups of natural numbers generated by three elements and its applications to algebraic geometry can be seen in \cite{Herzog70}. Some inexplicit formulas for the Frobenius number in three variables can be seen in \cite{tr17}.  

When $p>0$, the situation become even more difficult. For two variables, it is still easy to find explicit formulas of $g_p(a,b)$, $n_p(a,b)$ and $s_p(a,b)$. However, for three or more variables, no explicit formula had been found. But, finally in 2022 we succeeded in obtaining closed formulas for some special cases, including the triplets of triangular numbers \cite{Ko22a}, repunits \cite{Ko22b}, Fibonacci \cite{KY} and Jacobsthal numbers \cite{KP}.     

We are interested in finding any explicit closed formula for $p\ge 0$. In this paper, we give explicit formulas for the triples forming arithmetic progressions $A=\{a,a+d,a+2 d\}$, where $a$ and $d$ are positive integers with $\gcd(a,d)=1$. The main result is given as follows (Theorem \ref{th:3varp}). 
For $0\le p\le\fl{a/2}$, 
\begin{align*}
&g_p(a,a+d,a+2 d)=(a+2 d)p+\fl{\frac{a-2}{2}}a+(a-1)d\,,\\   
&n_p(a,a+d,a+2 d)\\
&=\begin{cases}
(2 a+2 d-1-p)p+\frac{(a-1)(a+2 d-1)}{4}&\text{if $a$ is odd};\\ 
(2 a+2 d-1-p)p+\frac{(a-1)(a+2 d-1)+1}{4}&\text{if $a$ is even}\,.
\end{cases} 
\end{align*} 
We also show their explicit closed formulas for the power sum $s^{(\mu)}_p(a,a+d,a+2 d)$ (Theorem \ref{th:spsum}) and the weighted sum $s_{\lambda,p}^{(\mu)}(a,a+d,a+2 d)$, defined by 
$$
s_{\lambda,p}^{(\mu)}(a_1,\dots,a_k):=\sum_{n\in\mathbb N_0\backslash S_p(a_1,\dots,a_k)}\lambda^n n^\mu
$$ 
(Theorem \ref{th:swpsum}).

\section{Preliminaries}  

We introduce an extension of the Ap\'ery set (see \cite{Apery}) in order to obtain the formulas for $g_p(A)$, $n_p(A)$ and $s_p(A)$. Without loss of generality, we assume that $a_1=\min(A)$. 

\begin{Def}  
Let $p$ be a nonnegative integer. For a set of positive integers $A=\{a_1,a_2,\dots,a_k\}$ with $\gcd(A)=1$ and $a_1=\min(A)$ we denote by 
$$
{\rm Ap}_p(A)={\rm Ap}_p(a_1,a_2,\dots,a_k)=\{m_0^{(p)},m_1^{(p)},\dots,m_{a_1-1}^{(p)}\}\,, 
$$ 
the $p$-Ap\'ery set of $A$, where each $m_i^{(p)}$ ($1\le i\le a_1$) satisfied the conditions: 
$$
{\rm (i)}\, m_i^{(p)}\equiv i\pmod{a_1},\quad{\rm (ii)}\, m_i^{(p)}\in S_p(A),\quad{\rm (iii)}\, m_i^{(p)}-a_1\not\in S_p(A)\,. 
$$ 
Note that $m_0^{(0)}$ is defined to be $0$.  
\label{apery} 
\end{Def}  

\noindent 
It follows that for each $p$, the set ${\rm Ap}_p(A)$ is a complete residue system modulo $a_1$. That is,  
$$
{\rm Ap}_p(A)\equiv\{0,1,\dots,a_1-1\}\pmod{a_1}\,. 
$$

By using the elements in ${\rm Ap}_p(A)$, the power sum of the elements in $\mathbb N_0\backslash S_p(A)$ can be given \cite{Ko-p} (see also \cite{Ko22z}).  
   
\begin{Prop}
Let $k$, $p$ and $\mu$ be integers with $k\ge 2$, $p\ge 0$ and $\mu\ge 0$.  
Assume that $\gcd(A)=1$.  We have 
\begin{align*} 
s^{(\mu)}_p(A)&:=\sum_{n\in\mathbb N_0\backslash S_p(A)}n^\mu\\ 
&=\frac{1}{\mu+1}\sum_{\kappa=0}^{\mu}\binom{\mu+1}{\kappa}B_{\kappa}a_1^{\kappa-1}\sum_{i=0}^{a_1-1}\bigl(m_i^{(p)}\bigr)^{\mu+1-\kappa}\\ 
&\qquad +\frac{B_{\mu+1}}{\mu+1}(a_1^{\mu+1}-1)\,,
\end{align*} 
where $B_n$ are Bernoulli numbers defined by 
$$
\frac{x}{e^x-1}=\sum_{n=0}^\infty B_n\frac{x^n}{n!}\,. 
$$ 
\label{prop1}
\end{Prop} 

When $\mu=0,1$ in Proposition \ref{prop1}, together with $g_p(A)$ we have formulas for the $p$-Frobenius number, the $p$-Sylvester number and the $p$-Sylvester sum.  

\begin{Lem}  
Let $k$ and $p$ be integers with $k\ge 2$ and $p\ge 0$.  
Assume that $\gcd(a_1,a_2,\dots,a_k)=1$.  We have 
\begin{align}  
g_p(A)&=\left(\max_{0\le j\le a_1-1}m_j^{(p)}\right)-a_1\,, 
\label{mp-g}
\\  
n_p(A)&=\frac{1}{a_1}\sum_{j=0}^{a_1-1}m_j^{(p)}-\frac{a_1-1}{2}\,,  
\label{mp-n}
\\
s_p(A)&=\frac{1}{2 a_1}\sum_{j=0}^{a_1-1}\bigl(m_j^{(p)}\bigr)^2-\frac{1}{2}\sum_{j=0}^{a_1-1}m_j^{(p)}+\frac{a_1^2-1}{12}\,.
\label{mp-s}
\end{align} 
\label{cor-mp}
\end{Lem} 

\noindent 
{\it Remark.}  
When $p=0$, the formulas (\ref{mp-g}), (\ref{mp-n}) and (\ref{mp-s}) reduce to the formulas by Brauer and Shockley \cite{bs62}, Selmer \cite{se77}, and Tripathi \cite{tr08}, respectively: 
\begin{align*}   
g(A)&=\left(\max_{1\le j\le a_1-1}m_j\right)-a_1\,,\\  
n(A)&=\frac{1}{a_1}\sum_{j=0}^{a_1-1}m_j-\frac{a_1-1}{2}\,,\\  
s(A)&=\frac{1}{2 a_1}\sum_{j=0}^{a_1-1}(m_j)^2-\frac{1}{2}\sum_{j=0}^{a_1-1}m_j+\frac{a_1^2-1}{12}\,,   
\end{align*} 
where $m_j=m_j^{(0)}$ ($1\le j\le a_1-1$) with $m_0=m_0^{(0)}=0$.

\section{The main result}

In this section, we shall show the main result and give its proof.  

\begin{theorem}  
Let $a$ and $d$ be integers with $a\ge 3$, $d>0$ and $\gcd(a,d)=1$.  
Then for $0\le p\le\fl{a/2}$, 
\begin{align*}
&g_p(a,a+d,a+2 d)=(a+2 d)p+\fl{\frac{a-2}{2}}a+(a-1)d\,,\\   
&n_p(a,a+d,a+2 d)\\
&=\begin{cases}
(2 a+2 d-1-p)p+\frac{(a-1)(a+2 d-1)}{4}&\text{if $a$ is odd};\\ 
(2 a+2 d-1-p)p+\frac{(a-1)(a+2 d-1)+1}{4}&\text{if $a$ is even}\,.
\end{cases} 
\end{align*} 
\label{th:3varp}
\end{theorem}  

\noindent 
{\it Remark.}  
When $p=0$, the formulas reduce to 
$$
g_0(a,a+d,a+2 d)=\fl{\frac{a-2}{2}}a+(a-1)d  
$$ 
and 
$$
n_0(a,a+d,a+2 d)=\begin{cases}
\frac{(a-1)(a+2 d-1)}{4}&\text{if $a$ is odd};\\ 
\frac{(a-1)(a+2 d-1)+1}{4}&\text{if $a$ is even}\,.
\end{cases}
$$ 
respectively, 
which are \cite[(3.9),(3.10)]{se77} when $k=3$.

Let $r_{x_2,x_3}=(a+d)x_2+(a+2 d)x_3$ for nonnegative integers $x_2$ and $x_3$. In the following tables, denote it by $(x_2,x_3)$ for simplicity.  

When $a$ is odd, as seen in \cite[(3.6)]{se77}, ${\rm Ap}_0(A_3)$ ($A_3=\{a,a+d,a+2d\}$) as the complete residue system (minimal system) modulo $a$ is given by Table \ref{tb:g00system}.  

\begin{table}[htbp]
  \centering
\begin{tabular}{cc}
\cline{1-2}
\multicolumn{1}{|c}{$(0,0)$}&\multicolumn{1}{c|}{$(1,0)$}\\
\multicolumn{1}{|c}{$(0,1)$}&\multicolumn{1}{c|}{$(1,1)$}\\
\multicolumn{1}{|c}{$\vdots$}&\multicolumn{1}{c|}{$\vdots$}\\
\multicolumn{1}{|c}{$(0,\frac{a-3}{2})$}&\multicolumn{1}{c|}{$(1,\frac{a-3}{2})$}\\
\cline{2-2}
\multicolumn{1}{|c|}{$(0,\frac{a-1}{2})$}&\\
\cline{1-1}
\end{tabular}
  \caption{Complete residue system ${\rm Ap}_0(A_3)$ for odd $a$}
  \label{tb:g00system}
\end{table} 

Concerning the complete residue system ${\rm Ap}_1(A_3)$, each congruent value modulo $a$ moves up one line to the upper right block. However, only the two values in the top row move to fill the gap below the first block. 
Namely, for $x_2\ge 1$
$$
r_{0,x_2}\equiv r_{2,x_2-1},\quad r_{1,x_2}\equiv r_{3,x_2-1}\pmod{a}
$$ 
and 
$$
r_{0,0}\equiv r_{1,\frac{a-1}{2}},\quad r_{1,0}\equiv r_{0,\frac{a+1}{2}}\pmod{a}\,.
$$ 
Since 
\begin{align*} 
3(a+d)+(a+2 d)x_3&=a+(a+d)+(a+2 d)(x_3+1)\,,\\
2(a+d)+(a+2 d)x_3&=a+(a+2 d)(x_3+1)\quad (x_3\ge 0)\,,\\ 
(a+d)+\frac{a-1}{2}(a+2 d)&=\left(\frac{a+1}{2}+d\right)a\,,\\ 
\frac{a+1}{2}(a+2 d)&=\left(\frac{a-1}{2}+d\right)a+(a+d)\,, 
\end{align*} 
we can know that each element in ${\rm Ap}_1(A_3)$ has exactly two expressions in terms of $(a,a+d,a+2 d)$. Note that each element minus $a$ has only one expression which is yielded from the right-hand side, that is, 
\begin{align*} 
&(a+d)+(a+2 d)(x_3+1),\quad (a+2 d)(x_3+1)\quad (x_3\ge 0),\\
&\left(\frac{a-1}{2}+d\right)a,\quad \left(\frac{a-3}{2}+d\right)a+(a+d)\,.  
\end{align*} 
Considering the maximal value, it is clear that 
$$
r_{0,\frac{a+1}{2}}>r_{1,\frac{a-1}{2}}>r_{2,\frac{a-3}{2}}>r_{3,\frac{a-5}{2}}\,.
$$

\begin{table}[htbp]
  \centering
\begin{tabular}{cccc}
\cline{1-2}\cline{3-4}
\multicolumn{1}{|c}{$(0,0)$}&\multicolumn{1}{c|}{$(1,0)$}&$(2,0)$&\multicolumn{1}{c|}{$(3,0)$}\\
\multicolumn{1}{|c}{$(0,1)$}&\multicolumn{1}{c|}{$(1,1)$}&$(2,1)$&\multicolumn{1}{c|}{$(3,1)$}\\
\multicolumn{1}{|c}{$\vdots$}&$\vdots$&\multicolumn{1}{|c}{$\vdots$}&\multicolumn{1}{c|}{$\vdots$}\\
\multicolumn{1}{|c}{$(0,\frac{a-5}{2})$}&\multicolumn{1}{c|}{$(1,\frac{a-5}{2})$}&$(2,\frac{a-5}{2})$&\multicolumn{1}{c|}{$(3,\frac{a-5}{2})$}\\
\cline{4-4}
\multicolumn{1}{|c}{$(0,\frac{a-3}{2})$}&\multicolumn{1}{c|}{$(1,\frac{a-3}{2})$}&\multicolumn{1}{c|}{$(2,\frac{a-3}{2})$}&\\
\cline{2-3}
\multicolumn{1}{|c|}{$(0,\frac{a-1}{2})$}&\multicolumn{1}{c|}{$(1,\frac{a-1}{2})$}&&\\
\cline{1-2}
\multicolumn{1}{|c|}{$(0,\frac{a+1}{2})$}&&&\\
\cline{1-1}
\end{tabular}
  \caption{Complete residue system ${\rm Ap}_1(A_3)$ from ${\rm Ap}_0(A_3)$ for odd $a$}
  \label{tb:g0system}
\end{table}

Concerning the complete residue system ${\rm Ap}_2(A_3)$, each congruent value modulo $a$ moves up one line to the upper right block (the third block). However, only the two values in the top row in the second block move to fill the gap below the first block. 
Namely, for $x_2\ge 1$
$$
r_{2,x_2}\equiv r_{4,x_2-1},\quad r_{3,x_2}\equiv r_{5,x_2-1}\pmod{a}
$$ 
and 
$$
r_{2,0}\equiv r_{1,\frac{a+1}{2}},\quad r_{3,0}\equiv r_{0,\frac{a+3}{2}}\pmod{a}\,.
$$ 
Since 
\begin{align*} 
5(a+d)+(a+2 d)x_3&=a+3(a+d)+(a+2 d)(x_3+1)\\ 
&=2 a+(a+d)+(a+2 d)(x_3+2)\,,\\
4(a+d)+(a+2 d)x_3&=a+2(a+d)+(a+2 d)(x_3+1)\\
&=2 a+(a+2 d)(x_3+2)\quad (x_3\ge 0)\,,\\ 
3(a+d)+\frac{a-3}{2}(a+2 d)&=a+(a+d)+\frac{a-1}{2}(a+2 d)\\
&=\left(\frac{a+3}{2}+d\right)a\,,\\
2(a+d)+\frac{a-1}{2}(a+2 d)&=a+\frac{a+1}{2}(a+2 d)\\
&=\left(\frac{a+1}{2}+d\right)a+(a+d)\,,\\ 
(a+d)+\frac{a+1}{2}(a+2 d)&=\left(\frac{a-1}{2}+d\right)a+2(a+d)\\
&=\left(\frac{a+1}{2}+d\right)a+(a+2 d)\,,\\
\frac{a+3}{2}(a+2 d)&=\left(\frac{a-3}{2}+d\right)a+3(a+d)\\
&=\left(\frac{a-1}{2}+d\right)a+(a+d)+(a+2 d)\,, 
\end{align*} 
we can know that each element in ${\rm Ap}_1(A_3)$ has exactly three expressions in terms of $(a,a+d,a+2 d)$. Note that each element minus $a$ has two expressions which are yielded from the right-hand side, that is, 
\begin{align*} 
3(a+d)+(a+2 d)(x_3+1)&=a+(a+d)+(a+2 d)(x_3+2)\,,\\
2(a+d)+(a+2 d)(x_3+1)&=a+(a+2 d)(x_3+2)\quad (x_3\ge 0)\,,\\
(a+d)+\frac{a-1}{2}(a+2 d)&=\left(\frac{a+1}{2}+d\right)a\,,\\
\frac{a+1}{2}(a+2 d)&=\left(\frac{a-1}{2}+d\right)a+(a+d)\,,\\
\left(\frac{a-3}{2}+d\right)a+2(a+d)&=\left(\frac{a-1}{2}+d\right)a+(a+2 d)\,,\\ 
\left(\frac{a-3}{2}+d-1\right)a+3(a+d)&=\left(\frac{a-3}{2}+d\right)a+(a+d)+(a+2 d)\,.  
\end{align*}   
Considering the maximal value, it is clear that 
$$
r_{0,\frac{a+3}{2}}>r_{1,\frac{a+1}{2}}>r_{2,\frac{a-1}{2}}>r_{3,\frac{a-3}{2}}>r_{4,\frac{a-5}{2}}>r_{5,\frac{a-7}{2}}\,. 
$$

\begin{table}[htbp]
  \centering
\begin{tabular}{cccccc}
\cline{1-2}\cline{3-4}\cline{5-6}
\multicolumn{1}{|c}{}&\multicolumn{1}{c|}{}&$(2,0)$&\multicolumn{1}{c|}{$(3,0)$}&$(4,0)$&\multicolumn{1}{c|}{$(5,0)$}\\
\multicolumn{1}{|c}{}&\multicolumn{1}{c|}{}&$(2,1)$&\multicolumn{1}{c|}{$(3,1)$}&$(4,1)$&\multicolumn{1}{c|}{$(5,1)$}\\
\multicolumn{1}{|c}{}&&\multicolumn{1}{|c}{$\vdots$}&\multicolumn{1}{c|}{$\vdots$}&$\vdots$&\multicolumn{1}{c|}{$\vdots$}\\
\multicolumn{1}{|c}{}&&\multicolumn{1}{|c}{$\vdots$}&\multicolumn{1}{c|}{$\vdots$}&$(4,\frac{a-7}{2})$&\multicolumn{1}{c|}{$(5,\frac{a-7}{2})$}\\
\cline{6-6}
\multicolumn{1}{|c}{}&\multicolumn{1}{c|}{}&$(2,\frac{a-5}{2})$&\multicolumn{1}{c|}{$(3,\frac{a-5}{2})$}&\multicolumn{1}{c|}{$(4,\frac{a-5}{2})$}&\\
\cline{4-5}
\multicolumn{1}{|c}{}&\multicolumn{1}{c|}{}&\multicolumn{1}{c|}{$(2,\frac{a-3}{2})$}&\multicolumn{1}{c|}{$(3,\frac{a-3}{2})$}&&\\
\cline{2-3}\cline{4-4}
\multicolumn{1}{|c|}{}&\multicolumn{1}{c|}{$(1,\frac{a-1}{2})$}&\multicolumn{1}{c|}{$(2,\frac{a-1}{2})$}&&&\\
\cline{1-2}\cline{3-3}
\multicolumn{1}{|c|}{$(0,\frac{a+1}{2})$}&\multicolumn{1}{c|}{$(1,\frac{a+1}{2})$}&&&&\\
\cline{1-2}
\multicolumn{1}{|c|}{$(0,\frac{a+3}{2})$}&&&&&\\
\cline{1-1}
\end{tabular}
  \caption{Complete residue system ${\rm Ap}_2(A_3)$ from ${\rm Ap}_1(A_3)$ for odd $a$}
  \label{tb:g0system}
\end{table}

If this process is continued for $p=3,4,\dots$, when $p=(a-1)/2$, the state shown in Table \ref{tb:gpsystem-odd} is reached. Here, the shaded cell parts show the elements of ${\rm Ap}_{\frac{a-1}{2}}(A_3)$. Elements with the same residues modulo $a$ move in the following positions according to $p=0,1,2,\dots,(a-3)/2,(a-1)/2$.  
\begin{align*}  
&(0,0)\rightarrow (1,\frac{a-1}{2})\rightarrow (3,\frac{a-3}{2})\rightarrow\cdots\rightarrow (a-4,2)\rightarrow (a-2,1)\\ 
&(1,0)\rightarrow (0,\frac{a+1}{2})\rightarrow (2,\frac{a-1}{2})\rightarrow\cdots\rightarrow (a-5,3)\rightarrow (a-3,2)\\ 
&(0,1)\rightarrow (2,0)\rightarrow (1,\frac{a+1}{2})\rightarrow\cdots\rightarrow (a-6,4)\rightarrow (a-4,3)\\ 
&(1,1)\rightarrow (3,0)\rightarrow (0,\frac{a+3}{2})\rightarrow\cdots\rightarrow (a-7,5)\rightarrow (a-5,4)\\ 
&\qquad\cdots\\
&(0,\frac{a-3}{2})\rightarrow (2,\frac{a-5}{2})\rightarrow (4,\frac{a-7}{2})\rightarrow\cdots\rightarrow (a-3,0)\rightarrow (1,a-2)\\ 
&(1,\frac{a-3}{2})\rightarrow (3,\frac{a-5}{2})\rightarrow (5,\frac{a-7}{2})\rightarrow\cdots\rightarrow (a-2,0)\rightarrow (0,a-1)\\ 
&(0,\frac{a-1}{2})\rightarrow (2,\frac{a-3}{2})\rightarrow (4,\frac{a-5}{2})\rightarrow\cdots\rightarrow (a-3,1)\rightarrow (a-1,0)
\end{align*}

\begin{table}[htbp]
  \centering
\scalebox{0.7}{
\begin{tabular}{ccccccccc}
\cline{1-2}\cline{3-4}\cline{7-8}\cline{9-9}
\multicolumn{1}{|c}{$(0,0)$}&\multicolumn{1}{c|}{$(1,0)$}&$(2,0)$&\multicolumn{1}{c|}{$(3,0)$}&$\cdots$&$\cdots$&\multicolumn{1}{|c}{$(a-3,0)$}&\multicolumn{1}{c|}{$(a-2,0)$}&\multicolumn{1}{|c|}{\cellcolor[gray]{0.8}$(a-1,0)$}\\
\cline{8-9}
\multicolumn{1}{|c}{$(0,1)$}&\multicolumn{1}{c|}{$(1,1)$}&$(2,1)$&\multicolumn{1}{c|}{$(3,1)$}&&&\multicolumn{1}{|c}{$(a-3,1)$}&\multicolumn{1}{|c|}{\cellcolor[gray]{0.8}$(a-2,1)$}&\\
\cline{7-8}
\multicolumn{1}{|c}{$\vdots$}&\multicolumn{1}{c|}{$\vdots$}&$\vdots$&\multicolumn{1}{c|}{$\vdots$}&&&\multicolumn{1}{|c|}{\cellcolor[gray]{0.8}$(a-3,2)$}&&\\
\cline{7-7}
\multicolumn{1}{|c}{$\vdots$}&\multicolumn{1}{c|}{$\vdots$}&$\vdots$&\multicolumn{1}{c|}{$\vdots$}&&&&&\\
\cline{4-4}
\multicolumn{1}{|c}{$(0,\frac{a-3}{2})$}&\multicolumn{1}{c|}{$(1,\frac{a-3}{2})$}&$(2,\frac{a-3}{2})$&\multicolumn{1}{|c|}{$(3,\frac{a-3}{2})$}&&&&&\\
\cline{2-2}\cline{3-4}
\multicolumn{1}{|c}{$(0,\frac{a-1}{2})$}&\multicolumn{1}{|c|}{$(1,\frac{a-1}{2})$}&\multicolumn{1}{c|}{$(2,\frac{a-1}{2})$}&&&&&&\\
\cline{1-2}\cline{3-3}
\multicolumn{1}{|c}{$(0,\frac{a+1}{2})$}&\multicolumn{1}{|c|}{$(1,\frac{a+1}{2})$}&&&&&&&\\
\cline{1-2}\cline{4-4}
\multicolumn{1}{|c|}{$(0,\frac{a+3}{2})$}&&&\multicolumn{1}{|c|}{\cellcolor[gray]{0.8}$(3,a-4)$}&&&&&\\
\cline{1-1}\cline{3-4}
$\vdots$&&\multicolumn{1}{|c|}{\cellcolor[gray]{0.8}$(2,a-3)$}&&&&&&\\
\cline{2-3}
$\vdots$&\multicolumn{1}{|c|}{\cellcolor[gray]{0.8}$(1,a-2)$}&&&&&&&\\
\cline{1-2}
\multicolumn{1}{|c|}{\cellcolor[gray]{0.8}$(0,a-1)$}&&&&&&&&\\
\cline{1-1}
\end{tabular}
} 
  \caption{Complete residue system ${\rm Ap}_p(A_3)$ for odd $a$}
  \label{tb:gpsystem-odd}
\end{table}

Indeed, each element $m_j^{(\frac{a-1}{2})}$ in ${\rm Ap}_{\frac{a-1}{2}}(A_3)$ has exactly $(a+1)/2$ expressions because   
\begin{align*}  
(a-1)(a+2 d)&=(d+i)a+(a-2 i)(a+d)+(i-1)(a+2 d)\\
&\qquad(1\le i\le\frac{a-1}{2})\,,\\
(a+d)+(a-2)(a+2 d)&=(d+i+1)a+(a-2 i-1)(a+d)\\ 
&\quad +(i-1)(a+2 d)\quad(1\le i\le\frac{a-1}{2})\,,\\
2(a+d)+(a-3)(a+2 d)&=a+(a-2)(a+2 d)\\
&=(d+i+2)a+(a-2 i-2)(a+d)\\
&\quad +(i-1)(a+2 d)\quad(1\le i\le\frac{a-3}{2})\,,\\
3(a+d)+(a-4)(a+2 d)&=a+(a+d)+(a-3)(a+2 d)\\
&=(d+i+3)a+(a-2 i-3)(a+d)\\
&\quad +(i-1)(a+2 d)\quad(1\le i\le\frac{a-3}{2})\,,\\
\cdots\\
(a-3)(a+d)+2(a+2 d)&=i a+(a-2i-3)(a+d)+(i+2)(a+2 d)\\
&\qquad(1\le i\le\frac{a-3}{2})\\
&=(a+d-2)a+(a+d)\,,\\
(a-2)(a+d)+(a+2 d)&=i a+(a-2i-2)(a+d)+(i+1)(a+2 d)\\
&\qquad(1\le i\le\frac{a-3}{2})\\
&=(a+d-1)a\,,\\ 
(a-1)(a+d)&=i a+(a-2i-1)(a+d)+i(a+2 d)\\
&\qquad(1\le i\le\frac{a-1}{2})\,.  
\end{align*} 
But any element $m_j^{(\frac{a-1}{2})}-a$ has $(a-1)/2$ expressions. 

Considering the maximal value, it is clear that 
$$
r_{0,a-1}>r_{1,a-2}>r_{2,a-3}>r_{3,a-4}>\cdots>r_{a-3,2}>r_{a-2,1}>r_{a-1,0}\,. 
$$ 

However, such a process until $0\le p\le(a-1)/2$ cannot be continued further.  When $p=(a+1)/2$, the same residue of $r_{a-1,0}$ modulo $a$ comes to the position $(1,a-1)$ and no element comes to the bottom left position $(0,a)$. Thus, after that, the pattern shifts and is complicated, so it becomes difficult to determine where the maximum element of ${\rm Ap}_p(A_3)$ is.  

Therefore, when $a$ is odd, for $0\le p\le(a-1)/2$, 
\begin{align*}
g_p(a,a+d,a+2 d)&=\left(\frac{a-1}{2}+p\right)(a+2 d)-a\\
&=(a+2 d)p+\frac{a(a-3)+2(a-1)d}{2}\,. 
\end{align*} 

Concerning the number of the representations, we need the summation of the elements in ${\rm Ap}_p(A_3)$. 
The elements in the staircase are  
\begin{multline*} 
\left(0,\frac{a-1}{2}+p\right),~\left(1,\frac{a-1}{2}+p-1\right),~\dots,~\\
\left(2 p-2,\frac{a-1}{2}-p+2\right),~\left(2 p-1,\frac{a-1}{2}-p+1\right)\,, 
\end{multline*}  
and the elements in the large block are  
$$
(2 p,0),~\dots,~\left(2 p,\frac{a-1}{2}-p\right),~(2 p+1,0),~\dots,~\left(2 p+1,\frac{a-1}{2}-p-1\right)\,. 
$$  
Hence, 
\begin{align}  
&\sum_{j=0}^{a-1}m_j^{(p)}\notag\\
&=\bigl(0+1+\cdots+(2 p-1)\bigr)(a+d)\notag\\
&\quad +\left(\left(\frac{a-1}{2}-p+1\right)+\left(\frac{a-1}{2}-p+2\right)+\cdots+\left(\frac{a-1}{2}+p\right)\right)(a+2 d)\notag\\
&\quad +2 p\left(\frac{a+1}{2}-p\right)(a+d)+\left(0+1+\cdots+\left(\frac{a-1}{2}-p\right)\right)(a+2 d)\notag\\
&\quad +(2 p+1)\left(\frac{a+1}{2}-p-1\right)(a+d)\notag\\
&\quad +\left(0+1+\cdots+\left(\frac{a-1}{2}-p-1\right)\right)(a+2 d)\notag\\
&=\frac{(2 p-1)(2 p)}{2}(a+d)\notag\\
&\quad +\left(\frac{(a-1+2 p)(a+1+2 p)}{8}-\frac{(a-1-2 p)(a+1-2 p)}{8}\right)(a+2 d)\notag\\
&\quad +p(a+1-2 p)(a+d)+\frac{(a-1-2 p)(a+1-2 p)}{8}(a+2 d)\notag\\
&\quad +\frac{2 p+1}{2}(a-2 p-1)(a+d)+\frac{(a-3-2 p)(a-1-2 p)}{8}(a+2 d)\notag\\
&=\frac{a}{4}\bigl((a+d)^2-(d+1)^2-4 p^2+4(2 a+2 d-1)p\bigr)\,. 
\label{eq:nodd} 
\end{align}
Hence, by Lemma \ref{cor-mp} (\ref{mp-n}) 
\begin{align*}
&n_p(a,a+d,a+2 d)\\
&=\frac{1}{4}\bigl((a+d)^2-(d+1)^2-4 p^2+4(2 a+2 d-1)p\bigr)-\frac{a-1}{2}\\
&=\frac{(a-1)(a+2 d-1)}{4}-p^2+(2 a+2 d-1)p\,.
\end{align*}

When $a$ is even, ${\rm Ap}_0(A_3)$ ($A_3=\{a,a+d,a+2d\}$) is given by the following figure.  

\begin{table}[htbp]
  \centering
\begin{tabular}{cc}
\cline{1-2}
\multicolumn{1}{|c}{$(0,0)$}&\multicolumn{1}{c|}{$(1,0)$}\\
\multicolumn{1}{|c}{$(0,1)$}&\multicolumn{1}{c|}{$(1,1)$}\\
\multicolumn{1}{|c}{$\vdots$}&\multicolumn{1}{c|}{$\vdots$}\\
\multicolumn{1}{|c}{$(0,\frac{a}{2}-1)$}&\multicolumn{1}{c|}{$(1,\frac{a}{2}-1)$}\\
\cline{1-2}
\end{tabular}
  \caption{Complete residue system ${\rm Ap}_0(A_3)$ for even $a$}
  \label{tb:g0system}
\end{table} 

Each congruent value modulo $a$ moves up one line to the upper right block. However, only the two values in the top row move to fill the gap below the first block. 
Namely, for $x_2\ge 1$
$$
r_{0,x_2}\equiv r_{2,x_2-1},\quad r_{1,x_2}\equiv r_{3,x_2-1}\pmod{a}
$$ 
and 
$$
r_{0,0}\equiv r_{0,\frac{a}{2}},\quad r_{1,0}\equiv r_{1,\frac{a}{2}}\pmod{a}\,.
$$ 

When the process continued for $p=1,2,\dots,a/2$, the state shown in Table \ref{tb:gpsystem-even} is reached. Here, the shaded cell parts show the elements of ${\rm Ap}_{\frac{a}{2}}(A_3)$. Elements with the same residues modulo $a$ move in the following positions according to $p=0,1,2,\dots,a/2-1,a/2$.  
\begin{align*}  
&(0,0)\rightarrow (1,\frac{a}{2})\rightarrow (2,\frac{a}{2}-1)\rightarrow\cdots\rightarrow (a-4,2)\rightarrow (a-2,1)\\ 
&(1,0)\rightarrow (1,\frac{a}{2})\rightarrow (3,\frac{a}{2}-1)\rightarrow\cdots\rightarrow (a-3,2)\rightarrow (a-1,1)\\ 
&(0,1)\rightarrow (2,0)\rightarrow (0,\frac{a}{2}+1)\rightarrow\cdots\rightarrow (a-6,4)\rightarrow (a-4,3)\\ 
&(1,1)\rightarrow (3,0)\rightarrow (1,\frac{a}{2}+1)\rightarrow\cdots\rightarrow (a-5,4)\rightarrow (a-3,3)\\ 
&\qquad\cdots\\
&(0,\frac{a}{2}-1)\rightarrow (2,\frac{a}{2}-2)\rightarrow (4,\frac{a}{2}-3)\rightarrow\cdots\rightarrow (a-2,0)\rightarrow (0,a-1)\\ 
&(1,\frac{a}{2}-1)\rightarrow (3,\frac{a}{2}-2)\rightarrow (5,\frac{a}{2}-3)\rightarrow\cdots\rightarrow (a-1,0)\rightarrow (1,a-1)
\end{align*}

\begin{table}[htbp]
  \centering
\scalebox{0.5}{
\begin{tabular}{cccccccccccc}
\cline{1-2}\cline{3-4}\cline{5-6}\cline{9-10}\cline{11-12}
\multicolumn{1}{|c}{$(0,0)$}&$(1,0)$&\multicolumn{1}{|c}{$(2,0)$}&$(3,0)$&\multicolumn{1}{|c}{$(4,0)$}&\multicolumn{1}{c|}{$(5,0)$}&$\cdots$&$\cdots$&\multicolumn{1}{|c}{$(a-4,0)$}&$(a-3,0)$&\multicolumn{1}{|c}{$(a-2,0)$}&\multicolumn{1}{c|}{$(a-1,0)$}\\ 
\cline{11-12} 
\multicolumn{1}{|c}{$(0,1)$}&$(1,1)$&\multicolumn{1}{|c}{$(2,1)$}&$(3,1)$&\multicolumn{1}{|c}{$(4,1)$}&\multicolumn{1}{c|}{$(5,1)$}&&&\multicolumn{1}{|c}{$(a-4,1)$}&$(a-3,1)$&\multicolumn{1}{|c}{\cellcolor[gray]{0.8}$(a-2,1)$}&\multicolumn{1}{c|}{\cellcolor[gray]{0.8}$(a-1,1)$}\\ 
\cline{9-10}\cline{11-12}
\multicolumn{1}{|c}{$\vdots$}&&\multicolumn{1}{|c}{$\vdots$}&\multicolumn{1}{c|}{}&$\vdots$&\multicolumn{1}{c|}{}&&&\multicolumn{1}{|c}{$(a-4,2)$}&\multicolumn{1}{c|}{$(a-3,2)$}&&\\  
\cline{9-10}
\multicolumn{1}{|c}{$\vdots$}&&\multicolumn{1}{|c}{$(2,\frac{a}{2}-3)$}&\multicolumn{1}{c|}{$(3,\frac{a}{2}-3)$}&$(4,\frac{a}{2}-3)$&\multicolumn{1}{c|}{$(5,\frac{a}{2}-3)$}&&&\multicolumn{1}{|c}{\cellcolor[gray]{0.8}$(a-4,3)$}&\multicolumn{1}{c|}{\cellcolor[gray]{0.8}$(a-3,3)$}&&\\ 
\cline{5-6}\cline{9-10}
\multicolumn{1}{|c}{$(0,\frac{a}{2}-2)$}&\multicolumn{1}{c|}{$(1,\frac{a}{2}-2)$}&$(2,\frac{a}{2}-2)$&\multicolumn{1}{c|}{$(3,\frac{a}{2}-2)$}&&&&&&&&\\ 
\cline{3-4}
\multicolumn{1}{|c}{$(0,\frac{a}{2}-1)$}&\multicolumn{1}{c|}{$(1,\frac{a}{2}-1)$}&$(2,\frac{a}{2}-1)$&\multicolumn{1}{c|}{$(3,\frac{a}{2}-1)$}&&&&&&&&\\ 
\cline{1-2}\cline{3-4} 
\multicolumn{1}{|c}{$(0,\frac{a}{2})$}&\multicolumn{1}{c|}{$(1,\frac{a}{2})$}&&&&&&&&&&\\ 
\cline{1-2}
\multicolumn{1}{|c}{$(0,\frac{a}{2}+1)$}&\multicolumn{1}{c|}{$(1,\frac{a}{2}+1)$}&&&&&&&&&&\\ 
\cline{1-2}
$\vdots$&&&&&&&&&&&\\ 
\cline{3-4}\cline{5-6}
$\vdots$&&\multicolumn{1}{|c}{$(2,a-5)$}&\multicolumn{1}{c|}{$(3,a-5)$}&\cellcolor[gray]{0.8}$(4,a-5)$&\multicolumn{1}{c|}{\cellcolor[gray]{0.8}$(5,a-5)$}&&&&&&\\
\cline{3-4}\cline{5-6}  
$\vdots$&&\multicolumn{1}{|c}{$(2,a-4)$}&\multicolumn{1}{c|}{$(3,a-4)$}&&&&&&&&\\ 
\cline{1-2}\cline{3-4}
\multicolumn{1}{|c}{$(0,a-3)$}&\multicolumn{1}{c|}{$(1,a-3)$}&\cellcolor[gray]{0.8}$(2,a-3)$&\multicolumn{1}{c|}{\cellcolor[gray]{0.8}$(3,a-3)$}&&&&&&&&\\ 
\cline{1-2}\cline{3-4}
\multicolumn{1}{|c}{$(0,a-2)$}&\multicolumn{1}{c|}{$(1,a-2)$}&&&&&&&&&&\\ 
\cline{1-2}  
\multicolumn{1}{|c}{\cellcolor[gray]{0.8}$(0,a-1)$}&\multicolumn{1}{c|}{\cellcolor[gray]{0.8}$(1,a-1)$}&&&&&&&&&&\\
\cline{1-2}  
\end{tabular}
} 
  \caption{Complete residue system ${\rm Ap}_{\frac{a}{2}}(A_3)$ for even $a$}
  \label{tb:gpsystem-even}
\end{table}

Each element $m_j^{(\frac{a}{2})}$ in ${\rm Ap}_{\frac{a}{2}}(A_3)$ has exactly $a/2+1$ expressions because   
\begin{align*}  
(a-1)(a+2 d)&=(d+i)a+(a-2 i)(a+d)+(i-1)(a+2 d)\\
&\qquad(1\le i\le\frac{a}{2})\,,\\
(a+d)+(a-1)(a+2 d)&=(d+i)a+(a-2 i+1)(a+d)\\
&\quad +(i-1)(a+2 d)\quad(1\le i\le\frac{a}{2})\,,\\
2(a+d)+(a-3)(a+2 d)&=a+(a-2)(a+2 d)\\
&=(d+i+2)a+(a-2 i-2)(a+d)\\
&\quad +(i-1)(a+2 d)\quad(1\le i\le\frac{a}{2}-1)\,,\\
3(a+d)+(a-3)(a+2 d)&=a+(a+d)+(a-2)(a+2 d)\\
&=(d+i+2)a+(a-2 i-1)(a+d)\\
&\quad +(i-1)(a+2 d)\quad(1\le i\le\frac{a}{2}-1)\,,\\
4(a+d)+(a-5)(a+2 d)&=a+2(a+d)+(a-4)(a+2 d)\\
&=2 a+(a-3)(a+2 d)\\
&=(d+i+4)a+(a-2 i-4)(a+d)\\
&\quad +(i-1)(a+2 d)\quad(1\le i\le\frac{a}{2}-2)\,,\\
5(a+d)+(a-5)(a+2 d)&=a+3(a+d)+(a-4)(a+2 d)\\
&=2 a+(a+d)+(a-3)(a+2 d)\\
&=(d+i+4)a+(a-2 i-3)(a+d)\\
&\quad +(i-1)(a+2 d)\quad(1\le i\le\frac{a}{2}-2)\,,\\
\cdots\\
(a-2)(a+d)+(a+2 d)&=i a+(a-2i-2)(a+d)+(i+1)(a+2 d)\\
&\qquad(1\le i\le\frac{a}{2}-1)\\
&=(a+d-1)a\,,\\ 
(a-1)(a+d)+(a+2 d)&=i a+(a-2i-1)(a+d)+(i+1)(a+2 d)\\
&\qquad(1\le i\le\frac{a-1}{2})\\
&=(a+d-1)a+(a+d)\,.  
\end{align*} 
But any element $m_j^{(\frac{a}{2})}-a$ has $a/2$ expressions. 

Since the maximal value in ${\rm Ap}_p(A_3)$ is at the position $(1,a/2-1+p)$, when $a$ is even, for $0\le p\le a/2$, 
\begin{align*}
g_p(a,a+d,a+2 d)&=(a+d)+\left(\frac{a}{2}-1+p\right)(a+2 d)-a\\
&=(a+2 d)p+\frac{a(a-2)+2(a-1)d}{2}\,. 
\end{align*}  
However, when $p=a/2+1$, no element comes to the position $(0,a)$ or $(1,a)$ because there is no element of ${\rm Ap}_{\frac{a}{2}}(A_3)$ in the top row. Hence, after $p=a/2+1$, the pattern is shifted and the situation becomes irregular and complicated.

Concerning the number of the representations, the elements in the staircase are  
\begin{align*} 
&\left(0,\frac{a}{2}+p-1\right),~\left(1,\frac{a}{2}+p-1\right),\\ 
&\left(2,\frac{a}{2}+p-3\right),~\left(3,\frac{a}{2}+p-3\right),\\ 
&\qquad\cdots\\
&\left(2 p-4,\frac{a}{2}-p+3\right),~\left(2 p-3,\frac{a}{2}-p+3\right),\\ 
&\left(2 p-2,\frac{a}{2}-p+1\right),~\left(2 p-1,\frac{a}{2}-p+1\right)\,, 
\end{align*}  
and the elements in the large block are  
\begin{align*} 
&(2 p,0),~(2 p,1),~\dots,~\left(2 p,\frac{a}{2}-p-1\right),\\
&(2 p+1,0),~(2 p+1,1),~\dots,~\left(2 p+1,\frac{a}{2}-p-1\right)\,. 
\end{align*}  
Hence, 
\begin{align}  
&\sum_{j=0}^{a-1}m_j^{(p)}\notag\\
&=\bigl(0+2+\cdots+(2 p-2)\bigr)(a+d)+\bigl(1+3+\cdots+(2 p-1)\bigr)(a+d)\notag\\
&\quad +2\left(\left(\frac{a}{2}-p+1\right)+\left(\frac{a}{2}-p+3\right)+\cdots+\left(\frac{a}{2}+p-1\right)\right)(a+2 d)\notag\\ 
&\quad +\left(\frac{a}{2}-p\right)(2 p)(a+d)+\left(\frac{a}{2}-p\right)(2 p+1)(a+d)\notag\\
&\quad +2\left(0+1+\cdots+\left(\frac{a}{2}-p-1\right)\right)(a+2 d)\notag\\ 
&=(2 p-1)p(a+d)+a p(a+2 d)\notag\\ 
&\quad +\left(\frac{a}{2}-p\right)(4 p+1)(a+d)+\left(\frac{a}{2}-p-1\right)\left(\frac{a}{2}-p\right)(a+2 d)\notag\\
&=\frac{a}{4}\bigl((a+d)^2-(d+1)^2+1-4 p^2+4(2 a+2 d-1)p\bigr)\,. 
\label{eq:neven}
\end{align}
Hence, by Lemma \ref{cor-mp} (\ref{mp-n}) 
\begin{align*}
&n_p(a,a+d,a+2 d)\\
&=\frac{1}{4}\bigl((a+d)^2-(d+1)^2+1-4 p^2+4(2 a+2 d-1)p\bigr)-\frac{a-1}{2}\\
&=\frac{(a-1)(a+2 d-1)+1}{4}-p^2+(2 a+2 d-1)p\,.
\end{align*}

In \cite[(3,9)]{se77}, 
$$
g\bigl(a,a+d,\dots,a+(k-1)d\bigr)=\fl{\frac{a-2}{k-1}}a+(a-1)d\,, 
$$ 
In \cite[(3,10)]{se77}, 
$$
n\bigl(a,a+d,\dots,a+(k-1)d\bigr)=\frac{(a-1)(q+d)+r(q+1)}{2}\,, 
$$ 
where integers $q$ and $r$ are determined as 
$$
a-1=q(k-1)+r,\quad 0\le r<k-1\,. 
$$

\section{Power sums}  

More generally, we can show a formula for the $p$-Sylvester power sum 
$$
s^{(\mu)}_p(a_1,a_2,\dots,a_k)=\sum_{d(n)\le p}n^\mu\quad(\mu\ge 1)\,, 
$$ 
so that $s_p(a_1,a_2,\dots,a_k)=s^{(1)}_p(a_1,a_2,\dots,a_k)$ and $n_p(a_1,a_2,\dots,a_k)$\\$=s^{(0)}_p(a_1,a_2,\dots,a_k)$.  
Once we know the exact structure of every element in ${\it Ap}_p(A)$, by applying Proposition \ref{prop1}, we can get the formula. Namely, we need to calculate $(m_i^{(p)}\bigr)^\nu$ for $\nu\ge 0$. From the previous section, when $n$ is odd, 
\begin{align*}  
(m_i^{(p)}\bigr)^\nu&=\sum_{j=0}^\nu\binom{\nu}{j}\sum_{k=0}^{2 p-1}k^{\nu-j}\left(\frac{a-1}{2}+p-k\right)^j(a+d)^{\nu-j}(a+2 d)^j\\
&\quad +\sum_{j=0}^\nu\binom{\nu}{j}(2 p)^{\nu-j}\sum_{k=0}^{\frac{a-1}{2}-p}k^j(a+d)^{\nu-j}(a+2 d)^j\\
&\quad +\sum_{j=0}^\nu\binom{\nu}{j}(2 p+1)^{\nu-j}\sum_{k=0}^{\frac{a-1}{2}-p-1}k^j(a+d)^{\nu-j}(a+2 d)^j\\
&=\sum_{j=0}^\nu\binom{\nu}{j}\Biggl(\sum_{k=0}^{2 p-1}k^{\nu-j}\left(\frac{a-1}{2}+p-k\right)^j\\
&\qquad +(2 p)^{\nu-j}\sum_{k=0}^{\frac{a-1}{2}-p}k^j+(2 p+1)^{\nu-j}\sum_{k=0}^{\frac{a-1}{2}-p-1}k^j
\Biggr)(a+d)^{\nu-j}(a+2 d)^j\,. 
\end{align*}
When $n$ is even 
\begin{align*}  
(m_i^{(p)}\bigr)^\nu&=\sum_{j=0}^\nu\binom{\nu}{j}\sum_{k=1}^{p}\bigl((2 k-2)^{\nu-j}+(2 k-1)^{\nu-j}\bigr)\left(\frac{a}{2}+p-2 k+1\right)^j\\
&\qquad\qquad\times(a+d)^{\nu-j}(a+2 d)^j\\
&\quad +\sum_{j=0}^\nu\binom{\nu}{j}(2 p)^{\nu-j}\sum_{k=0}^{\frac{a}{2}-p-1}k^j(a+d)^{\nu-j}(a+2 d)^j\\
&\quad +\sum_{j=0}^\nu\binom{\nu}{j}(2 p+1)^{\nu-j}\sum_{k=0}^{\frac{a}{2}-p-1}k^j(a+d)^{\nu-j}(a+2 d)^j\\
&=\sum_{j=0}^\nu\binom{\nu}{j}\Biggl(\sum_{k=1}^{p}\bigl((2 k-2)^{\nu-j}+(2 k-1)^{\nu-j}\bigr)\left(\frac{a}{2}+p-2 k+1\right)^j\\
&\qquad +\bigl((2 p)^{\nu-j}+(2 p+1)^{\nu-j}\bigr)\sum_{k=0}^{\frac{a}{2}-p-1}k^j\Biggr)(a+d)^{\nu-j}(a+2 d)^j\,. 
\end{align*}
Hence, by Proposition \ref{prop1}, we obtain the generalized Sylvester power sum for $(a,a+d,a+2 d)$. 

\begin{theorem}  
Let $a$, $d$, $p$ and $\mu$ be integers with $a\ge 3$, $d>0$, $p\ge 0$, $\mu\ge 1$ and $\gcd(a,d)=1$.  Then, when $a$ is odd, we have 
\begin{align*} 
&s^{(\mu)}_p(a,a+d,a+2 d)\\ 
&=\frac{1}{\mu+1}\sum_{\kappa=0}^{\mu}\binom{\mu+1}{\kappa}B_{\kappa}a^{\kappa-1}\sum_{i=0}^{a-1}\sum_{j=0}^{\mu+1-\kappa}\binom{\mu+1-\kappa}{j}\\
&\qquad\times\Biggl(\sum_{k=0}^{2 p-1}k^{\mu+1-\kappa-j}\left(\frac{a-1}{2}+p-k\right)^j\\
&\qquad\quad +(2 p)^{\mu+1-\kappa-j}\sum_{k=0}^{\frac{a-1}{2}-p}k^j+(2 p+1)^{\mu+1-\kappa-j}\sum_{k=0}^{\frac{a-1}{2}-p-1}k^j
\Biggr)\\
&\qquad\qquad\times(a+d)^{\mu+1-\kappa-j}(a+2 d)^j\\  
&\quad +\frac{B_{\mu+1}}{\mu+1}(a^{\mu+1}-1)\,.  
\end{align*} 
When $a$ is even, we have 
\begin{align*} 
&s^{(\mu)}_p(a,a+d,a+2 d)\\ 
&=\frac{1}{\mu+1}\sum_{\kappa=0}^{\mu}\binom{\mu+1}{\kappa}B_{\kappa}a^{\kappa-1}\sum_{i=0}^{a-1}\sum_{j=0}^{\mu+1-\kappa}\binom{\mu+1-\kappa}{j}\\
&\qquad\times\Biggl(\sum_{k=1}^{p}\bigl((2 k-2)^{\mu+1-\kappa-j}+(2 k-1)^{\mu+1-\kappa-j}\bigr)\left(\frac{a}{2}+p-2 k+1\right)^j\\
&\qquad\quad +\bigl((2 p)^{\mu+1-\kappa-j}+(2 p+1)^{\mu+1-\kappa-j}\bigr)\sum_{k=0}^{\frac{a}{2}-p-1}k^j\Biggr)\\
&\qquad\qquad\times(a+d)^{\mu+1-\kappa-j}(a+2 d)^j\\ 
&\quad+\frac{B_{\mu+1}}{\mu+1}(a^{\mu+1}-1)\,.  
\end{align*}   
\label{th:spsum} 
\end{theorem}

In particular, when $\mu=2$, Theorem \ref{th:spsum} reduces the formula of the $p$-Sylvester sum of the triple forming an arithmetic progression.  

\begin{Cor}
Let $a$ and $d$ be integers with $a\ge 3$, $d>0$ and $\gcd(a,d)=1$.
Then for $0\le p\le\fl{a/2}$, when $a$ is odd, we have 
\begin{align*}  
s_p(a,a+d,a+2 d)
&=\frac{(a-1)(a+2 d-1)(a^2+2 a d-a-d-2)}{24}\\
&\quad +\frac{3 a^3+9 a^2(d-1)+2 a(3 d^2-9 d+1)-6 d^2+2 d}{6}p\\ 
&\quad +\frac{3 a^2+a(6 d-1)+4 d^2-d}{2}p^2-\frac{3(a+d)}{2}p^3\,. 
\end{align*}
When $a$ is even, we have 
\begin{align*}  
s_p(a,a+d,a+2 d)
&=\frac{(a-1)(a+2 d-1)(a^2+2 a d-a-d-2)+3(a^2+2 a d-d)}{24}\\
&\quad +\frac{3 a^3+9 a^2(d-1)+a(6 d^2-18 d+5)-6 d^2+5 d}{6}p\\ 
&\quad +\frac{3 a^2+a(6 d-1)+4 d^2-d}{2}p^2-\frac{3(a+d)}{2}p^3\,. 
\end{align*}
\label{cor:ssum}
\end{Cor}

\section{Weighted sums} 

In this section, we consider the weighted sums whose numbers of representations are less than or equal to $p$ \cite{KZ0,KZ}:  
$$
s_{\lambda,p}^{(\mu)}(a_1,\dots,a_k):=\sum_{n\in{\it Ap}_p(a_1,\dots,a_k)}\lambda^n n^\mu\,, 
$$ 
where $\lambda\ne 1$ and $\mu$ is a positive integer.   

Here,   
Eulerian numbers $\eul{n}{m}$ appear in the generating function 
\begin{equation}
\sum_{k=0}^\infty k^n x^k=\frac{1}{(1-x)^{n+1}}\sum_{m=0}^{n-1}\eul{n}{m}x^{m+1}\quad(n\ge 1)
\label{eu:gf}
\end{equation} 
with $0^0=1$ and $\eul{0}{0}=1$, 
and have an explicit formula 
$$
\eul{n}{m}=\sum_{k=0}^{m}(-1)^k\binom{n+1}{k}(m-k+1)^n\,. 
$$  
Then, an explicit formula of the $p$-weighted sum is given in terms of the elements in ${\it Ap}_p(a_1,\dots,a_k)$ \cite{Ko-p} (see also \cite{Ko22z}).   

\begin{Lem}  
Assume that $\lambda\ne 1$ and $\lambda^{a_1}\ne 1$. Then for a positive integer $\mu$,  
\begin{align*}  
&s_{\lambda,p}^{(\mu)}(a_1,\dots,a_k)\\
&=\sum_{n=0}^\mu\frac{(-a_1)^n}{(\lambda^{a_1}-1)^{n+1}}\binom{\mu}{n}\sum_{j=0}^n\eul{n}{n-j}\lambda^{j a_1}\sum_{i=0}^{a_1-1}\bigl(m_i^{(p)}\bigr)^{\mu-n}\lambda^{m_i^{(p)}}\\
&\quad +\frac{(-1)^{\mu+1}}{(\lambda-1)^{\mu+1}}\sum_{j=0}^\mu\eul{\mu}{\mu-j}\lambda^j\,.
\end{align*}
\label{lem-hh}
\end{Lem} 

In order to obtain the formula for the $p$-Sylvester weighted power sum, we need to calculate $(m_i^{(p)}\bigr)^\nu\lambda^{m_i^{(p)}}$ for $\nu\ge 0$, $\lambda^a\ne 1$ and $\lambda^d\ne 1$. We use the formula  
$$
a^\nu x^a=\sum_{i=0}^\nu\sts{\nu}{i}x^i\frac{d^i}{d x^i}x^a\,, 
$$ 
where $\sts{n}{m}$ denote the Stirling numbers of the second kind, calculated as 
$$
\sts{n}{m}=\frac{1}{m!}\sum_{i=0}^m(-1)^i\binom{m}{i}(m-i)^n\,. 
$$ 
From the previous section, when $n$ is odd, 
\begin{align*} 
&(m_i^{(p)}\bigr)^\nu\lambda^{m_i^{(p)}}\\
&=\sum_{i=0}^\nu\sts{\nu}{i}\lambda^i\left[\frac{d^i}{d x^i}\left(\sum_{l=0}^{2 p-1}x^{l(a+d)+(\frac{a-1}{2}+p-l)(a+2 d)}\right.\right.\\ 
&\quad\left.\left. +\sum_{l=0}^{\frac{a-1}{2}-p}x^{2 p(a+d)+l(a+2 d)}+\sum_{l=0}^{\frac{a-1}{2}-p-1}x^{(2 p+1)(a+d)+l(a+2 d)}\right)\right]_{x=\lambda}\\
&=\sum_{i=0}^\nu\sts{\nu}{i}\lambda^i\left[\frac{d^i}{d x^i}\left(\frac{x^{a(\frac{a-1}{2}+d+p)}(x^{2 d p}-1)}{x^d-1}\right.\right.\\ 
&\quad\left.\left. +\frac{x^{2 p(a+d)}(x^{(a+2 d)(\frac{a+1}{2}-p)}-1)}{x^{a+2 d}-1}+\frac{x^{(2 p+1)(a+d)}(x^{(a+2 d)(\frac{a-1}{2}-p)}-1)}{x^{a+2 d}-1}\right)\right]_{x=\lambda}\,.
\end{align*} 
When $n$ is even, 
\begin{align*} 
&(m_i^{(p)}\bigr)^\nu\lambda^{m_i^{(p)}}\\
&=\sum_{i=0}^\nu\sts{\nu}{i}\lambda^i\left[\frac{d^i}{d x^i}\left(\sum_{l=0}^{p-1}x^{2 l(a+d)+(\frac{a}{2}+p-2 l-1)(a+2 d)}\right.\right.\\ 
&\quad +\sum_{l=0}^{p-1}x^{(2 l+1)(a+d)+(\frac{a}{2}+p-2 l-1)(a+2 d)}\\ 
&\quad\left.\left. +\sum_{l=0}^{\frac{a}{2}-p-1}x^{2 p(a+d)+l(a+2 d)}+\sum_{l=0}^{\frac{a}{2}-p-1}x^{(2 p+1)(a+d)+l(a+2 d)}\right)\right]_{x=\lambda}\\ 
&=\sum_{i=0}^\nu\sts{\nu}{i}\lambda^i\left[\frac{d^i}{d x^i}\left(\frac{x^{\frac{a}{2}(a+2 d-2+2 p)}(x^{a+d}+1)(x^{2 d p}-1)}{x^{2 d}-1}\right.\right.\\ 
&\quad\left.\left. +\frac{x^{2(a+d)p}(x^{a+d}+1)(x^{(a+2 d)(\frac{a}{2}-p)}-1)}{x^{a+2 d}-1}\right)\right]_{x=\lambda}\,. 
\end{align*} 

Hence, by Lemma \ref{lem-hh}, we obtain the generalized Sylvester weighted power sum for $(a,a+d,a+2 d)$. 

\begin{theorem}  
Let $a$, $d$, $p$, $\lambda$ and $\mu$ be integers with $a\ge 3$, $d>0$, $p\ge 0$, $\lambda\ne 1$, $\lambda^d\ne 1$, $\mu\ge 1$ and $\gcd(a,d)=1$.  Then, when $a$ is odd, we have 
\begin{align*} 
&s_{\lambda,p}^{(\mu)}(a,a+d,a+2 d)\\
&=\sum_{n=0}^\mu\sum_{j=0}^n\sum_{i=0}^\nu\frac{(-a)^n}{(\lambda^{a}-1)^{n+1}}\binom{\mu}{n}\eul{n}{n-j}\sts{\nu}{i}\lambda^{j a+i}\\
&\quad\times\left[\frac{d^i}{d x^i}\left(\frac{x^{a(\frac{a-1}{2}+d+p)}(x^{2 d p}-1)}{x^d-1}\right.\right.\\ 
&\qquad\left.\left. +\frac{x^{2 p(a+d)}(x^{(a+2 d)(\frac{a+1}{2}-p)}-1)}{x^{a+2 d}-1}+\frac{x^{(2 p+1)(a+d)}(x^{(a+2 d)(\frac{a-1}{2}-p)}-1)}{x^{a+2 d}-1}\right)\right]_{x=\lambda}\\ 
&\quad +\frac{(-1)^{\mu+1}}{(\lambda-1)^{\mu+1}}\sum_{j=0}^\mu\eul{\mu}{\mu-j}\lambda^j\,.
\end{align*}  
When $a$ is even, we have 
\begin{align*} 
&s_{\lambda,p}^{(\mu)}(a,a+d,a+2 d)\\
&=\sum_{n=0}^\mu\sum_{j=0}^n\sum_{i=0}^\nu\frac{(-a)^n}{(\lambda^{a}-1)^{n+1}}\binom{\mu}{n}\eul{n}{n-j}\sts{\nu}{i}\lambda^{j a+i}\\
&\quad\times\left[\frac{d^i}{d x^i}\left(\frac{x^{\frac{a}{2}(a+2 d-2+2 p)}(x^{a+d}+1)(x^{2 d p}-1)}{x^{2 d}-1}\right.\right.\\ 
&\qquad\left.\left. +\frac{x^{2(a+d)p}(x^{a+d}+1)(x^{(a+2 d)(\frac{a}{2}-p)}-1)}{x^{a+2 d}-1}\right)\right]_{x=\lambda}\\ 
&\quad +\frac{(-1)^{\mu+1}}{(\lambda-1)^{\mu+1}}\sum_{j=0}^\mu\eul{\mu}{\mu-j}\lambda^j\,.
\end{align*}   
\label{th:swpsum} 
\end{theorem} 

\noindent 
{\it Remark.}  
The case $\lambda=1$ is not included in Theorem \ref{th:swpsum}, but in Theorem \ref{th:spsum}.

\section{Examples}  

For $(11,15,19)$, that is, $a=11$ and $d=4$, when $q=5$, by Theorem \ref{th:3varp}, Corollary \ref{cor:ssum}, Theorems \ref{th:spsum} and \ref{th:swpsum}, we have 
\begin{align*}  
g_{5}(11,15,19)&=179\,,\\ 
n_{5}(11,15,19)&=165\,,\\ 
s_{5}(11,15,19)&=13605\,,\\ 
s_{5}^{(3)}(11,15,19)&=189158535\quad(\mu=3)\,,\\ 
s_{2,5}^{(3)}(11,15,19)&=46691295420476497563538523364517263\\
&\qquad 55433630648909109181546522\quad(\lambda=2)\,. 
\end{align*}

For $(6,11,16)$, that is, $a=6$ and $d=5$, when $q=3$, we have 
\begin{align*}  
g_{3}(6,11,16)&=85\,,\\ 
n_{3}(6,11,16)&=73\,,\\ 
s_{3}(6,11,16)&=2675\,,\\ 
s_{3}^{(3)}(6,11,16)&=7652009\quad(\mu=3)\,,\\ 
s_{2,3}^{(3)}(6,11,16)&=24083450837052351738334815453210\quad(\lambda=2)\,. 
\end{align*}

In fact, the integers whose representations in terms of $(6,11,16)$ are less than or equal to $3$ are 
$$
0,\underbrace{1,2,\dots,59}_{59},61,62,63,64,65,67,68,69,73,74,75,79,85\,. 
$$  
Hence, for example, 
\begin{align*}  
s_{2,3}^{(3)}(6,11,16)&=2^0\cdot 0^3+2^1\cdot 1^3+2^2\cdot 2^3+\cdots+2^{59}\cdot 59^3+2^{61}\cdot 61^3\\
&\quad +2^{62}\cdot 62^3+2^{63}\cdot 63^3+2^{64}\cdot 64^3+2^{65}\cdot 65^3\\
&\quad +2^{67}\cdot 67^3+2^{68}\cdot 68^3+2^{69}\cdot 69^3+2^{73}\cdot 73^3\\
&\quad +2^{74}\cdot 74^3+2^{75}\cdot 75^3+2^{79}\cdot 79^3+2^{85}\cdot 85^3\\ 
&=24083450837052351738334815453210\,. 
\end{align*}

\section{Final comments}  

It should not be thought that the similar repetitive process to [KLP,KP,KY] is simply going on in any triple. For example, it is known that some triples of Pell sequences do not follow a similar process but the formation of the elements of $0$-Ap\'ery set is different (in preparation). 
 
In addition, it is still very difficult to find any explicit formula for four or more variables in the sequence of arithmetic progressions too. For the moment only in the case of repunits \cite{Ko22b}, for $p\ge 0$ explicit formulas about four and five repunits are obtained, though the structures are even more complicated than that of three variables. 
One of the reasons of the difficulties lies in the following. 
In the case of three variables, for any $j$ ($j=0,1,\dots,a-1$), $m_j^{(0)}<m_j^{(1)}<\dots$. However, in the case of four and more variables, for some $j$'s, $m_j^{(p)}=m_j^{(p+1)}$. It means that some elements in ${\rm Ap}_p(A)$ and in ${\rm Ap}_{p+1}(A)$ are overlapped.  

Selmer \cite{se77} found a formula of the Frobenius number for almost arithmetic sequences by generalizing the previous result (Roberts \cite{ro56} for $h=1$; Brauer \cite{br42} for $h=d=1$). For a positive integer $h$, 
$$
g_0(a,h a+d,\dots,h a+(k-1)d)=\left(h\fl{\frac{a-2}{k-1}}+h-1\right)a+(a-1)d\,.
$$ 
Selmer also gave an explicit formula for the Sylvester number $n_0(a,h a+d,\dots,h a+(k-1)d)$.    
Some formulas for the Sylvester sum $s_0(a,h a+d,\dots,h a+(k-1)d)$ and its variations are given in \cite{Ko22}.  
However, it is known that even when $d=2$ (the sequence of consecutive odd numbers), we have not found any explicit form of $g_p(a,a+2,\dots,h a+2(k-1))$ for general $p>0$.  

Another problem is whether we can find any convenient formula when $p>\fl{a/2}$.



\section*{Conflict of interests} 

There is no conflict of interests.



\end{document}